\documentclass[11pt]{article}

\usepackage{latexsym}
\usepackage{amssymb}

\newtheorem{theorem}{Theorem}

\newtheorem{corollary}{Corollary}

\begin{document}
{
\begin{center}
{\Large\bf
On bounded complex Jacobi matrices and related moment problems.}
\end{center}
\begin{center}
{\bf Sergey M. Zagorodnyuk}
\end{center}

\noindent
\textbf{Abstract.}
In this paper we study the linear functional $S$ on complex polynomials which is associated to a bounded complex Jacobi matrix $J$. 
The associated moment problem is considered: find a positive Borel measure $\mu$ on $\mathbb{C}$ subject to conditions
$\int z^n d\mu = s_n$, where $s_n$ are prescribed complex numbers (moments). This moment problem may be viewed as an
extension of the Stieltjes and Hamburger moment problems to the complex plane.
Sufficient conditions for the solvability of the moment problem are provided.
As a corollary, we obtain conditions 
for the existence of an integral representation $S(p) = \int_\mathbb{C} p(z) d\mu$, with a positive Borel measure $\mu$.
An interrelation of the associated to the complex Jacobi matrix operator $A_0$, acting in $l^2$ on finite vectors,
and the multiplication by z operator in $L^2_\mu$ is discussed as well. 



\section{Introduction.}

The theory of real Jacobi matrices is a well-known and classical subject with a lot of applications in various domains of
mathematics and other sciences, see the books of Akhiezer and Berezanskii~\cite{cit_2000_Akhiezer_book_1965},\cite{cit_7000_Berezanskii}.
Complex Jacobi matrices or J-matrices appeared in a context of J-fractions, see Wall's book~\cite{cit_2500_Wall}.
They have not attracted so much attention as their real versions. 
An important recent work on complex Jacobi matrices was done by Beckermann in 2001, who collected and arranged in a nice form basic facts
on this subject, see~\cite{cit_2100_Beckermann_2001}.
After Beckermann's paper the study of complex Jacobi matrices became essentially more active. We can mention the following
directions of investigations:
perturbations and spectral analysis (see~\cite{cit_2050_Swiderski_2020} and references therein);
quadrature rules (\cite{cit_2030_M_C_2003}); eigenvalue problems (\cite{cit_Ikebe_1996});
determinacy questions (\cite{cit_2150_Beckermann_Castro_Smirnova_2004}).
Two-sided Jacobi matrices are also studied intensively: for the real case we refer to~\cite{cit_7000_Berezanskii},
and for recent developments see~\cite{cit_2035_S_2017},\cite{cit4000_Ch_Z_2017} and references therein.

By a complex Jacobi matrix one means a semi-infinite tridiagonal complex matrix of the following form:

\begin{equation}
\label{f1_5}
J =
\left(
\begin{array}{ccccc}
b_0 & a_0 & 0 & 0 & \ldots \\
a_0 & b_1 & a_1 & 0 & \ldots \\
0 & a_1 & b_2 & a_2 & \ldots \\
\vdots & \vdots & \vdots & \vdots & \ddots 
\end{array}
\right),
\end{equation}
where $a_k,b_k\in\mathbb{C}:$ $a_k\not=0$, $k\in\mathbb{Z}_+$.
Let us recall some basic known facts about complex Jacobi matrices which we shall need in what follows.
By matrix multiplication the matrix $J$ generates a linear operator $A_0$ on $l^2_{fin}$. If
\begin{equation}
\label{f1_7}
|a_k| \leq M,\quad |b_k| \leq M,\qquad \forall k\in\mathbb{Z}_+,\quad  \mbox{for some $M>0$},
\end{equation}
then $A_0$ is bounded. In this case, by continuity it extends on the whole space $l_2$ to a bounded operator $A$~\cite{cit_2100_Beckermann_2001}.

With a complex Jacobi matrix $J$ one associates a system of polynomials $\{ p_n(\lambda) \}_{n=0}^\infty$, $p_0(\lambda) = 1$,
such that
\begin{equation}
\label{f1_10}
a_{n-1} p_{n-1}(\lambda) + b_n p_n(\lambda) + a_n p_{n+1}(\lambda) = \lambda p_n(\lambda),\qquad n\in\mathbb{Z}_+,
\end{equation}
where $a_{-1}:=0$, $p_{-1}(\lambda):=0$.
A linear with respect to the both arguments functional $\sigma(u,v)$, $u,v\in\mathbb{P}$, which satisfies
\begin{equation}
\label{f1_12}
\sigma(p_n(\lambda),p_m(\lambda)) = \delta_{n,m},\qquad  n,m\in\mathbb{Z}_+,
\end{equation}
is said to be \textit{the spectral function} of the difference equation
\begin{equation}
\label{f1_14}
a_{n-1} y_{n-1} + b_n y_n + a_n y_{n+1} = \lambda y_n,\qquad n\in\mathbb{Z}_+,
\end{equation}
see~\cite{cit_2500_Z_Serdica}. 
The difference equation~(\ref{f1_14}) (and therefore the complex Jacobi matrix $J$) can be recovered by its
spectral function and a sequence of signs, see~\cite{cit_2500_Z_Serdica} for the details of this reconstruction.

\begin{theorem}(\cite[Theorem 1]{cit_2500_Z_Serdica})
\label{t1_1}
A linear with respect to the both arguments functional $\sigma(u,v)$, $u,v\in\mathbb{P}$,
is the spectral function of a difference equation of type~(\ref{f1_14}) iff:

1) $\sigma(\lambda u(\lambda), v(\lambda)) = \sigma(u(\lambda), \lambda v(\lambda)),\qquad u,v\in\mathbb{P}$;

2) $\sigma(1,1) = 1$;

3) For arbitrary polynomial $u_k(\lambda)$ of degree $k$, there exists a polynomial $\widehat u_k(\lambda)$ of degree $k$ such that:
$$ \sigma(u_k(\lambda),\widehat u_k(\lambda)) \not= 0. $$
\end{theorem}

By property~1) we see that
\begin{equation}
\label{f1_16}
\sigma(u(\lambda),v(\lambda)) = \sigma(u(\lambda)v(\lambda), 1),\qquad u,v\in\mathbb{P}.
\end{equation}
Consider the following linear functional $S$, which is said to be \textit{associated to the complex Jacobi matrix $J$}:
\begin{equation}
\label{f1_18}
S(u) = \sigma(u, 1),\qquad u\in\mathbb{P}.
\end{equation}
By~(\ref{f1_12}) it has the following property:
\begin{equation}
\label{f1_20}
S(p_n p_m) = \delta_{n,m},\qquad  n,m\in\mathbb{Z}_+.
\end{equation}
Denote
\begin{equation}
\label{f1_22}
\mathbf{s}_n := S(\lambda^n),\qquad  n\in\mathbb{Z}_+.
\end{equation}
The numbers $\{ \mathbf{s}_n \}_{n=0}$ are said to be \textit{the moments} of $S$.
Our main objective here is to provide conditions on the moments $\mathbf{s}_n$, which imply the existence of an integral representation
of $S$ of the following form:
\begin{equation}
\label{f1_24}
S(p) = \int_\mathbb{C} p(z) d\mu,
\end{equation}
with a (non-negative) Borel measure $\mu$.
We shall use the following moment problem:
find a (non-negative) measure $\mu$ on $\mathfrak{B}(\mathbb{C})$
such that
\begin{equation}
\label{f1_50}
\int_{\mathbb{C}} z^k d\mu(z) = s_k,\qquad  k\in\mathbb{Z}_+.
\end{equation}
Here $\{ s_k \}_{k\in\mathbb{Z}_+}$ is a prescribed set of complex numbers (moments).
This moment problem was recently stated in~\cite{cit_3000_Zagorodnyuk_2022_Similarity}.
Sufficient conditions for the solvability of the moment problem~(\ref{f1_50}) will be given in
Theorem~\ref{t2_1}. In a consequence, some conditions 
for the existence of an integral representation~(\ref{f1_24}) with a positive Borel measure $\mu$
will appear in Corollary~\ref{c2_1}.
As another application, we have a representation of the operator $A_0$ as a multiplication by $z$ operator in 
$L^2_\mu$, see Corollary~2.

Finally, we remark that the operator $A$ is complex symmetric (see e.g.~\cite{cit5000} for definitions), and it belongs to 
the class $C_+(H)$, for $H = l^2$, see~\cite{cit_3500_Z}.

\noindent
{\bf Notations. }
As usual, we denote by $\mathbb{R}, \mathbb{C}, \mathbb{N}, \mathbb{Z}, \mathbb{Z}_+$
the sets of real numbers, complex numbers, positive integers, integers and non-negative integers,
respectively. 
By $\mathbb{Z}_{k,l}$ we mean all integers $r$, which satisfy the following inequality:
$k\leq r\leq l$.
By $\mathbb{P}$ we mean a set of all complex polynomials.
By $\mathfrak{B}(M)$ we denote the set of all Borel subsets of a set $M\subseteq\mathbb{C}$.
For a measure $\mu$ on $\mathfrak{B}(M)$ we denote by $L^2_\mu = L^2_\mu(M)$ the usual space of all
(classes of equivalence of) Borel measurable complex-valued functions $f$ on $M$, such that
$\int_M |f|^2 d\mu < +\infty$.
The class of the equivalence containing a function $f$ will be denoted by $[f]$.
By $l^2$ we denote the usual space of square-summable complex sequences $\vec u = (u_k)_{k=0}^\infty$, $u_k\in\mathbb{C}$,
and $l^2_{fin}$ means the subset of all finite vectors from $l^2$. Moreover $\vec e_k$ means a vector
from $l^2$ having $1$ in $k$'s place and zeros in other places ($k\in\mathbb{Z}_+$).

If H is a Hilbert space then $(\cdot,\cdot)_H$ and $\| \cdot \|_H$ mean
the scalar product and the norm in $H$, respectively. 
Indices may be omitted in obvious cases.
For a linear operator $A$ in $H$, we denote by $D(A)$
its  domain, by $R(A)$ its range, and $A^*$ means the adjoint operator
if it exists. If $A$ is invertible then $A^{-1}$ means its
inverse. $\overline{A}$ means the closure of the operator, if the
operator is closable. If $A$ is bounded then $\| A \|$ denotes its
norm.
By $E_H$ we denote the identity operator in $H$, i.e. $E_H x = x$,
$x\in H$. In obvious cases we may omit the index $H$. If $H_1$ is a subspace of $H$, then $P_{H_1} =
P_{H_1}^{H}$ denotes the orthogonal projection of $H$ onto $H_1$.
By $A|_{H_1}$ we mean the restriction of $A$ to the subspace $H_1$.

\section{Moment problems on $\mathbb{C}$ and complex Jacobi matrices.}

At first we shall study the moment problem~(\ref{f1_50}) which may be viewed as an extension
of the Stieltjes moment problem (SMP) and the Hamburger moment problem (HMP). While the extension 
$$ \mbox{SMP } \rightarrow \mbox{HMP} $$
is well known, the extension to the complex plane was usually accompanied by adding additional monomials under the integral sign and
the corresponding moments (the complex moment problem).

\begin{theorem}
\label{t2_1}
Let the moment problem~(\ref{f1_50}) be given with some complex moments $\{ s_k \}_{k\in\mathbb{Z}_+}$, $s_0=1$.
If the following condition holds:
\begin{equation}
\label{f2_5}
|s_n| \leq R^n,\qquad n\in\mathbb{Z}_+,
\end{equation}
for some $R>0$, then the moment problem~(\ref{f1_50}) is solvable.
Moreover, in this case it has a solution $\mu$ with a compact support.
\end{theorem}
\textbf{Proof.} 
At first we shall consider the moment problem~(\ref{f1_50}) with moments $\{ s_k \}_{k\in\mathbb{Z}_+}$, $s_0=1$, which
satisfy the following condition:
\begin{equation}
\label{f2_7}
\sum_{n=0}^\infty |s_n|^2 < \infty.
\end{equation}
Introduce the following vectors:
\begin{equation}
\label{f2_9}
\vec g := 
\left(
\begin{array}{ccccc} 0 \\
\overline{s_1} \\
\overline{s_2} \\
\overline{s_3} \\
\vdots \end{array}
\right),\quad
\vec f := 
\left(
\begin{array}{ccccc} \overline{s_1} \\
\overline{s_2} - \overline{s_1}^2 \\
\overline{s_3} - \overline{s_2}\overline{s_1} \\
\overline{s_4} - \overline{s_3}\overline{s_1} \\
\vdots \end{array}
\right) =
\left(
\begin{array}{ccccc} \overline{s_1} \\
\overline{s_2} \\
\overline{s_3} \\
\overline{s_4} \\
\vdots \end{array}
\right) -
\overline{s_1}
\left(
\begin{array}{ccccc} 0 \\
\overline{s_1} \\
\overline{s_2} \\
\overline{s_3} \\
\vdots \end{array}
\right).
\end{equation}
By condition~(\ref{f2_7}) it follows that $\vec f,\vec g\in l^2$.
Let $S$ be the right shift operator on $l^2$:
$$ S \vec u = (0,u_0,u_1,...),\qquad  \vec u = (u_0,u_1,...)\in l^2. $$
Consider the following operator $B$, which is defined on the whole $l^2$:
\begin{equation}
\label{f2_14}
B = S + (\cdot,\vec f)_{l^2} \vec e_0 - (\cdot,\vec g)_{l^2} \vec e_1.  
\end{equation}
The matrix $\mathcal{M}$ of the bounded operator $B$ with respect to the basis $\{ \vec e_k \}_{k=0}^\infty$ has the following form:
\begin{equation}
\label{f2_16}
\mathcal{M} = 
\left(
\begin{array}{ccccccc} 
s_1 & s_2 - s_1^2 & s_3 - s_2 s_1 & s_4 - s_3 s_1 & \cdots & s_{n+1} - s_n s_1 & \cdots \\
1 & -s_1 & -s_2 & - s_3 & \cdots & - s_n & \cdots \\
0 & 1 & 0 & 0 & \cdots & 0 & \cdots \\
0 & 0 & 1 & 0 & \cdots & 0 & \cdots \\
0 & 0 & 0 & 1 & \cdots & 0 & \cdots \\
\vdots & \vdots & \vdots & \vdots & \vdots & \vdots & \ddots \end{array}
\right).
\end{equation}
Consider the following vectors in $l^2$ (a similar construction was used in~\cite{cit_3540_Zagorodnyuk_2022_Axioms}):
\begin{equation}
\label{f2_18}
x_0 =  \vec e_0,\quad  x_j =  \vec e_j + s_j \vec e_0,\qquad j\in\mathbb{N},
\end{equation}
and the following operator on $l^2_{fin}$:
\begin{equation}
\label{f2_20}
M \sum_{k=0}^d \alpha_k x_k =  
\sum_{k=0}^d \alpha_k x_{k+1},\qquad \alpha_k\in\mathbb{C},\ d\in\mathbb{Z}_+.
\end{equation}
By induction we conclude that
\begin{equation}
\label{f2_22}
x_n = M^n x_0,\qquad n\in\mathbb{Z}_+,
\end{equation}
and therefore
\begin{equation}
\label{f2_24}
s_n = (x_n,x_0)_{l^2} = (M^n x_0,x_0)_{l^2},\qquad n\in\mathbb{Z}_+.
\end{equation}
Since
$$ \vec e_0 = x_0,\quad \vec e_j = x_j - s_j x_0,\qquad j\in\mathbb{N}, $$
then
$$ M \vec e_0 = \vec e_1 + s_1 \vec e_0, $$
$$ M \vec e_j = \vec e_{j+1} - s_j \vec e_1 + (s_{j+1} - s_j s_1) \vec e_0,\qquad j\in\mathbb{N}. $$
By a direct calculation one may verify that the matrix of $M$ with respect to the basis $\{ \vec e_k \}_{k=0}^\infty$
is exactly the matrix $\mathcal{M}$.
Therefore the operator $M$ is bounded and it extends on the whole space $l^2$ to the bounded operator $B$.
Denote $\rho := \| B \|$. Then the operator $\widetilde B := \frac{1}{\rho} B$ is a contraction. It has a unitary dilation $U$ in a Hilbert space
$\widehat H\supseteq l^2$ (see, e.g.,~\cite{cit_995_Sz.-Nagy_Book}):
\begin{equation}
\label{f2_26}
P^{\widehat H}_{l^2} U^k |_{l^2} = \widetilde B^k,\qquad k\in\mathbb{Z}_+.
\end{equation}
By~(\ref{f2_24}),(\ref{f2_26}) we may write:
\begin{equation}
\label{f2_28}
s_n = \rho^n ({\widetilde B}^n x_0, x_0)_{l^2} = \rho^n (U^n x_0, x_0)_{\widetilde H} = 
( (\rho U)^n x_0, x_0)_{\widetilde H},\ n\in\mathbb{Z}_+. 
\end{equation}
Since $\rho U$ is a bounded normal operator, its spectral resolution $E(\delta)$, $\delta\in\mathfrak{B}(\mathbb{C})$
provides a solution $\mu = (E(\delta) x_0,x_0)$ to the moment problem.

Consider now the moment problem from the assumptions of the theorem.
Choose an arbitrary $\tau > R$, and set
$$ \widetilde s_n := \frac{s_n}{ \tau^n },\qquad n\in\mathbb{Z}_+. $$
Then
$$ |s_n| \leq \left(  \frac{R}{\tau} \right)^n, $$
and therefore condition~(\ref{f2_7}) holds for $\widetilde s_n$.
Thus, by the already proved result we can see from~(\ref{f2_28}) that
$$ \widetilde s_n = ( N^n x_0, x_0)_{\widetilde H},\ n\in\mathbb{Z}_+, $$
for a bounded normal operator $N$ in a Hilbert space $\widetilde H\supseteq l^2$.
Therefore
$$ s_n = ( (\tau N)^n x_0, x_0)_{\widetilde H},\ n\in\mathbb{Z}_+. $$
The required result now follows from the spectral theorem for the bounded normal operator $\tau N$.
$\Box$

\begin{corollary}
\label{c2_1}
Let $J$ be a complex Jacobi matrix~(\ref{f1_5}). Suppose that condition~(\ref{f1_7}) holds. For the associated to $J$ linear functional $S$
the following integral representation holds:
\begin{equation}
\label{f2_30}
S(p) = \int_\mathbb{C} p(\lambda) d\mu,\qquad  p\in\mathbb{P}, 
\end{equation}
with a (non-negative) measure $\mu$ on $\mathfrak{B}(\mathbb{C})$, and $\mu$ has a compact support.
\end{corollary}
\textbf{Proof. }
In fact, for a complex Jacobi matrix one may write (see~\cite{cit_2100_Beckermann_2001}):
\begin{equation}
\label{f2_32}
\vec e_n = p_n(A) \vec e_0,\qquad n\in\mathbb{Z}_+,
\end{equation}
and therefore
\begin{equation}
\label{f2_34}
(p_m(A) \vec e_0, p_n(A) \vec e_0) = \delta_{m,n},\qquad m,n\in\mathbb{Z}_+,
\end{equation}
and
\begin{equation}
\label{f2_36}
S(u) = (u(A) \vec e_0, \vec e_0),\qquad u\in\mathbb{P}.
\end{equation}
Thus, for the moments $\mathbf{s}_n$ of $S$ we have the following estimate:
\begin{equation}
\label{f2_38}
|\mathbf{s}_n| \leq | (A^n \vec e_0, \vec e_0) | \leq \| A \|^n,\qquad n\in\mathbb{Z}_+.
\end{equation}
By Theorem~\ref{t2_1} it follows that there exists a compactly supported Borel measure $\mu$ on $\mathbb{C}$
such that
$$ S(z^n) = \mathbf{s}_n = \int z^n d\mu,\qquad  n\in\mathbb{Z}_+. $$
By linearity we obtain relation~(\ref{f2_30}). $\Box$

For a compactly supported measure $\mu$ on $\mathfrak{B}(\mathbb{C})$, we shall denote by $\Lambda_\mu$ the operator of multiplication
by the independent variable in $L^2_\mu$, and $\mathcal{P} := \{ [p],\ p\in\mathbb{P} \}$.

\begin{corollary}
\label{c2_2}
Let $J$ be a complex Jacobi matrix~(\ref{f1_5}) and $\{ p_n \}_{n=0}^\infty$, $p_0=1$,
be a system of polynomials satisfying~(\ref{f1_10}). Suppose that condition~(\ref{f1_7}) holds.
Let $A_0$ be the associated to $J$ operator on $l^2_{fin}$, and
\begin{equation}
\label{f2_40}
T \sum_{k=0}^d \xi_k \vec e_k = \left[ \sum_{k=0}^d \xi_k p_k(z) \right],\qquad \xi_k\in\mathbb{C},\ d\in\mathbb{Z}_+, 
\end{equation}
maps $l^2_{fin}$ into $L^2_\mu$, while $\mu$ is a positive measure provided by Corollary~\ref{c2_1}.
Then the linear operator $T$ is invertible and
\begin{equation}
\label{f2_42}
T A_0 T^{-1} = \Lambda_0,
\end{equation}
where $\Lambda_0=\Lambda_\mu|_\mathcal{P}$.
\end{corollary}
\textbf{Proof. }
Notice that
$$ T A_0 \vec e_k = T (b_{k-1} \vec e_{k-1} + a_k \vec e_k + b_k \vec e_{k+1}) = 
\left[ b_{k-1} p_{k-1}(z) + a_k p_k(z) + b_k p_{k+1}(z) \right] = $$
$$ = \left[ z p_k(z) \right] = \Lambda_0 [p_k] = \Lambda_0 T \vec e_k,\qquad k\in\mathbb{Z}_+. $$
By the linearity it follows that
\begin{equation}
\label{f2_44}
T A_0 = \Lambda_0 T.
\end{equation}
It remains to check that $T$ is invertible. Suppose to the contrary that there exists a nonzero vector 
$$ \vec u = \sum_{k=0}^r \eta_k \vec e_k \in l^2_{fin},\qquad   \eta_k\in\mathbb{C},\ \eta_r\not=0,\ r\in\mathbb{Z}_+,    $$ 
such that
$$ T \vec u = 0. $$
Then
$$ 0 = \| T \vec u \|^2 = \int \left| \sum_{k=0}^r \eta_k p_k(z) \right|^2 d\mu. $$
This means that the measure $\mu$ is finitely atomic, with atoms among the zeros of the polynomial
$u_r(z) := \sum_{k=0}^r \eta_k p_k(z)$, $\deg u_r=r$. 
By property~3) of Theorem~\ref{t1_1} there exists
a polynomial $\widehat u_r(\lambda)$ of degree $r$ such that:
$$ \int u_r(z) \widehat u_r(z) d\mu  = \sigma(u_r,\widehat u_r) \not= 0. $$
However the integral on the left side is equal to zero, a contradiction. 
The proof of the corollary is complete.
$\Box$

Of course, it would be valuable to obtain for the operator $A$ a result similar to that given in Corollary~\ref{c2_2} 
for the operator $A_0$.
However, it is not clear when $T$ is bounded and has a bounded inverse, except for the real matrices.
This question will be studied elsewhere.

}

\noindent
Address:

V. N. Karazin Kharkiv National University \newline\indent
School of Mathematics and Computer Sciences \newline\indent
Department of Higher Mathematics and Informatics \newline\indent
Svobody Square 4, 61022, Kharkiv, Ukraine

Sergey.M.Zagorodnyuk@gmail.com; zagorodnyuk@karazin.ua

\end{document}